\documentclass[11pt]{article}

\usepackage{comment}
\usepackage[T1]{fontenc}

\usepackage[utf8]{inputenc}

\usepackage[a4paper, total={6in, 8in}]{geometry}
\usepackage{array,graphicx}

\title{ CMALight: a novel Minibatch Algorithm for large-scale non convex finite sum optimization
}

\author{ Corrado Coppola \footnote{Department of Computer, Control and Management Engineering A. Ruberti - Sapienza, University of Rome. Via Ariosto 25 - 00185 Rome.} \footnote{Corresponding author. Email: corrado.coppola@uniroma1.it}, 
Giampaolo Liuzzi \footnotemark[1] , Laura Palagi \footnotemark[1] 
}

\usepackage{hyperref}

\usepackage{amsmath,amsfonts,amsthm,amssymb}	
\usepackage{todonotes}
\usepackage{bbm}
\usepackage{comment}
\usepackage{siunitx }
\usepackage{subfigure}
\usepackage{algpseudocode}
\usepackage{algorithm}
\usepackage{rotating}
\usepackage{adjustbox}
\usepackage{graphicx}
\usepackage{bigstrut}
\usepackage{multicol}
\usepackage{multirow}
\usepackage[english]{babel}
\usepackage[T1]{fontenc}
\usepackage{caption}
\usepackage{numprint}
\usepackage{footnote}

\usepackage[margins]{trackchanges}
\usepackage{todonotes}
\addeditor{RS}
\addeditor{LP}
\addeditor{GL}

\usepackage{todonotes}
\definecolor{oran}{rgb}{1.,0.4,0.}
\definecolor{dgreen}{rgb}{0.1,0.6,0.3}

\usepackage{nomencl}

\newtheorem{proposition}{Proposition}
\newtheorem{lemma}[proposition]{Lemma}

\newtheorem{assumption}{Assumption}

\newcommand{\R}{{\mathbb R}}

\newcommand{\CMALight}{\texttt{CMA Light}}

\renewcommand{\Re}{\mathbb{R}}

\newcommand{\commento}[1]{}
\newcommand{\CMA}{{\texttt{CMA }}}

\usepackage{curve2e}
\definecolor{gray}{gray}{0.4}

\begin{document}

\maketitle

\begin{abstract}
    The supervised training of a deep neural network on a given dataset consists in the unconstrained minimization of the finite sum of continuously differentiable functions, commonly referred to as loss with respect to the samples. These functions depend on the network parameters and most of the times are non-convex.  We develop \CMALight , a new globally convergent mini-batch gradient method to tackle this problem. We consider the recently introduced Controlled Minibatch Algorithm (CMA) framework and we overcome its main bottleneck, removing the need for at least one evaluation of the whole objective function per iteration. We prove global convergence of \CMALight under mild assumptions and we discuss extensive computational results on the same experimental test-bed used for CMA, showing that \CMALight requires less computational effort than most of the state-of-the-art optimizers. Eventually, we present early results on a large-scale Image Classification task.

\end{abstract}

\section{Introduction}
\label{sec:intro}
In this paper we consider the problem of the unconstrained minimization of the sum of P coercive, continuously differentiable functions $f_p$:
\begin{equation}
    \label{eq:minfw}
    \min_{w \in \R^n} f(w) = \sum_{i=1}^P f_i(w)
\end{equation}
We also consider that the functions $f_p$ may not be convex. While the problem itself is not particularly hard to describe from a theoretical point of view, computational issues arise when the number of functions P grows to indefinitely large values, which often occurs in many applications, such as constraints minimization violation \cite{mayne1981solving,zhang2010nonmonotone}, Stochastic Programming \cite{bollapragada2018adaptive,bottou2018optimization}, distributed optimization in Sensor Networks \cite{bertsekas2011incremental}.

Nonetheless, Problem \eqref{eq:minfw}
is mostly known since it represents the optimization problem behind the training of Machine Learning predictive models \cite{shalev2014understanding} on a given datasets of samples and labels $\{(x,y)_p, p=1,…,P\}$. 
In this case, it is required to solve the unconstrained minimization problem of a given loss function with respect to the network parameters. The most widespread loss functions are inherently expressed as a sum over P terms, e.g., the regularized Mean Square Error \eqref{eq:MSE}  for regression models and the regularized Cross Entropy Loss \eqref{eq:CEL} for classification models.
\begin{equation}
    \label{eq:MSE} f(w) = \frac 1 P \sum_{i=1}^P  (\hat y(w;x^i) - y^i)^2 + \lambda \vert \vert w \vert \vert^2
\end{equation}

\begin{equation}
    \label{eq:CEL}
    f(w) -\frac 1 P \sum_{i=1}^{P} [y^i \log( \hat y (x^i;w)) + (1 - y_i) \log(1 - \hat y (x^i;w))]+\lambda \|w\|^2 
\end{equation}

Training a Deep neural network on huge datasets results in a large-scale optimization problem, which makes computationally prohibitive to apply any standard first-order method (also referred as full-batch methods), requiring the evaluation of the whole $\nabla f(w)$ at each iteration. 
Computing $\nabla f(w)$ with respect to all variables and to all samples has a computational complexity of $o (nP)$, where $n$ can be in the order of trillions ($10^{12})$ for Vision Transformers \cite{khan2022transformers} and Diffusion Models \cite{croitoru2023diffusion}, while the number of samples $P$ can be in the order of terabytes \cite{lee2017big} in Big Data applications.
This is a computationally prohibitive challenge even for the most up-to-date devices.
Unfortunately, this issue also makes  impossible to rely on the strong convergence properties of these methods, which usually require only the Lipschitz-smoothness and coerciveness assumptions \cite{grippo2023basic}.

Therefore, Machine Learning community usually tackles Problem \eqref{eq:minfw} using methods that cycles over the samples and exploits only a small subset of samples $p^k$ (mini-batch) to perform a weights update. 
These methods are commonly referred to as mini-batch methods. 
They are classified according to the sampling rule of $p^k$. 
Following the historical development of these methods, we firstly recall the Robbins and Monroe 1952 SGD method \cite{robbins1951stochastic,Robbins&Monro:1951}, also known as Stochastic Gradient (SG), which picks up $p^k$ randomly. 
In Incremental Gradient (IG) \cite{bertsekas2011incremental} samples are picked up according to a pre-selected order, while in Random Reshuffling \cite{gurbuzbalaban2021random,shamir2016without} they are selected without replacement, according to an index permutation, changing every time the algorithm passes over all the P samples.

Despite being computationally very effective, MB methods convergence is very hard to prove and strong assumptions are usually required.

More in details, the most commonly used condition to prove SGD convergence is the growth condition on the second moment of the norm of the gradients with respect to all the samples \cite{bottou2018optimization}, i.e.,
$$ \mathbb{E}(\vert \vert \nabla f_i (w) \vert \vert ^2) \le a + b \vert \vert \nabla f (w) \vert \vert ^2, \ \ \forall i = 1, \dots,P
$$
Despite there have been many attemps in the literature to relax and to study more in details the convergence properties of SGD (see, e.g., \cite{schmidt2017minimizing,NIPS2014_5258,gower2021sgd,gower2021stochastic,lei2019stochastic,chen2018lag,duchi2018stochastic,asi2019stochastic,nguyen2022finite,Bottou10large-scalemachine,ghadimi2013stochastic}), either growth conditions or convexity assumptions are required, which is very hard to guarantee in Machine Learning training problems.

Concerning IG, convergence can be proved on a very strong growth condition \cite{bertsekas2011incremental}, i.e.,
$$ \| \nabla f_i (w) \| \le a + b \| \nabla f(w) \|, \ \ \forall i = 1, \dots,P
$$
which is verified only for specific classes of functions.

Eventually, the convergence theory of RR is the less developed, since the idea of permutation-based sampling rule introduces dependencies among the sampled gradients that are almost impossible to analyze from a theoretical point of view without introducing further assumptions, such as strong convexity of $f(w)$  \cite{gurbuzbalaban2021random,shamir2016without}.

However, despite the lack of convergence theory, the most used and computationally effective state-of-the-art algorithms for Machine Learning training are the Adaptive Gradient methods, which are based on adjusting the sampled gradients with the estimates computed in previous iterations, such as Adam and Adamax \cite{Kingma2015AdamAM}, Adagrad \cite{Duchi2010AdaptiveSM,Duchi:2011:ASM:1953048.2021068}, Adadelta \cite{Zeiler2012ADADELTAAA}, and NAdam \cite{Dozat2016IncorporatingNM}.

Controlled Minibatch Algorithm (CMA) is an algorithmic framework developed in \cite{liuzzi2022convergence} seeking for a trade-off balance between the efficiency standards required by the challenge of Big Data and the strong convergence properties ensured by full-batch method.

CMA is an ease-controlled modification of RR, where the global convergence is proved under mild assumptions (only coerciveness of $f(w)$ and Lipschitz-smoothness of $\nabla f_i(w)$). After each RR epoch, the algorithm checks through a watchdog condition, whether a sufficient decrease of the objective function has been attained, i.e., whether the tentative direction provided by the RR epoch is a "sufficiently good" descent direction.
When the condition does not hold, the step-size is either increased, either decreased, or driven to zero to restart the epoch, according to further conditions, involving the norm of the direction.
CMA has been proved to be computationally highly effective, especially on large-scale regression datasets and to our knowledge is the only globally convergent Minibatch gradient method under such mild assumptions.
Nonetheless, checking the sufficient decrease watchdog condition requires at least one evaluation of $f(w)$ per each epoch, which can still be computationally prohibitive for the most challenging tasks in Machine Learning, such as Image Classification and Video Recognition, where both the number of samples, and the number of features can indefinitely explode, as well as the network dimension.

The main objective of this paper is to overcome the CMA bottleneck, developing a new algorithm, which we will refer to as \CMALight, exploiting an estimate of $f(w)$ inherently provided by any MB gradient method, i.e. the sum of the batch losses, to check the sufficient decrease condition.
We firstly introduce our new algorithm and, under proper assumptions made in Section \ref{sec:assumptions}, we prove the global convergence of \CMALight in Section \ref{sec:cmalight}. In Section 4 we describe our computational experience in testing \CMALight against CMA and some early results on larger datasets and deeper networks that cannot be managed by CMA.
Finally, in Section 5 we present our preliminary conclusions concerning our research.

\section{Initial assumptions and background}
\label{sec:assumptions}
In this paper, we focus on ease-controlled modifications of Incremental Gradient (IG) and Random Reshuffling (RR) algorithms. Following the \CMA framework introduced in \cite{liuzzi2022convergence}, we aim at improving its efficiency making only seldomly necessary to evaluate the objective function on the whole dataset. 
While this is a common practice in Deep Learning community, especially when dealing with large-scale datasets and ultra-deep architectures, to the best of our knowledge no attempts to prove global converge have been made yet. 
We introduce the following assumptions on the objective function and on its components $f_p$.

\begin{assumption}\label{ass:compact}
The function $f(w)$ is coercive, i.e. all the level sets 
\[
  {\cal L}(w_0) = \{w\in\Re^n:\ f(w) \leq f(w_0)\}
\]
are compact for any $w_0\in\Re^n$.
\end{assumption}

\begin{assumption}\label{ass:lipschtizagradient}
The gradients $\nabla f_p(w)$ are Lipschitz continuous for each $p=1,\dots, P$, namely there exists $L>0$ such that
\[
 \|\nabla f_p(u)-\nabla f_p(v)\|\le L\|u-v\|\qquad \forall \ u,w\in\R^n
\]
\end{assumption}

\begin{assumption} \label{ass:fpcoercive}
The components $f_p$ are all non-negative and coercive, for each $p=1,\dots, P$, namely,
\begin{enumerate}
    \item $f_p(w) \ge 0$ for any $w \in \R^n$
    \item ${\cal L}(w_0) = \{w\in\Re^n:\ f_p(w) \leq f_p(w_0)\}$ is a compact set for each $p=1,\dots, P$ and for any $w_0 \in \R^n$
\end{enumerate}
\end{assumption}

Despite we add the assumption \ref{ass:fpcoercive}, which was not necessary in \cite{liuzzi2022convergence}, we remark that this holds at least for the most common loss functions used to train deep neural networks, e.g., the mean square error and the cross-entropy loss.

Similarly to \cite{liuzzi2022convergence}, our algorithm is based on a mini-batch gradient method iteration embedded in an outer iteration (epoch)).

\begin{algorithm}[htb]
\begin{algorithmic}[1]
    \State Given the sequences $\{\zeta^k\}$, and $\{w^k\}$ with $\zeta^k \in \R_+$ and $w^k\in\R^n$
    \For {$k=0,1,2\dots $}       
    \State Compute $(\widetilde w^k, d^k)$ = \texttt{Inner\_Cycle}($w^k,\zeta^k$)
    \EndFor
\end{algorithmic}
\caption{mini-batch Gradient (MG) iteration}
\label{alg:IG}
\end{algorithm}

Algorithm \ref{alg:IG} encompasses not only the Incremental Gradient method \cite{Bertsekas2000,blatt2007convergent,bertsekas2011incremental} but also most of the state-of-the-art optimizers used to train deep neural networks, such as  Adam \cite{Kingma2015AdamAM}, NAdam \cite{Dozat2016IncorporatingNM}, SGD \cite{robbins1951stochastic, bottou2018optimization, ruder2016overview, Sutskever2013OnTI}, Adagrad \cite{Duchi2010AdaptiveSM}, and Adadelta \cite{Zeiler2012ADADELTAAA}. However, due to its robust convergence properties \cite{Bertsekas2000}, we focus on the Incremental Gradient method. 
Since IG requires strong properties on the gradient of the objective functions to be satisfied, \CMA framework \cite{liuzzi2022convergence} has been proposed as an ease-controlled modification of IG method and proved to converge to a stationary point. 
Nonetheless, the watchdog rule introduced in the \CMA framework to force the sequence generated by the algorithm to remain within a compact level set requires at least one evaluation of the objective function per each epoch to verify that $f(w^k) \le f(w^0)$. While computational effort stays acceptable when training simple neural architectures on the regression and classification dataset composing the \CMA test-bed in \cite{liuzzi2022convergence}, evaluating the objective function on the whole dataset is computationally unaffordable for large-scale datasets and ultra-deep architectures. 

The underlying idea of \CMALight is to replace the real objective function with an estimation without introducing any further assumption on its behaviour. For this purpose we define the \texttt{Inner\_Cycle} in Algorithm \ref{alg:innercycle1}, which is by definition an epoch of IG method returning an approximated value of the objective function $\tilde f^k$.

\begin{algorithm}[h!]
\begin{algorithmic}[1]
    \State {\bf Input}: $w^k$, $\zeta^k$
    \State Set $\widetilde w_0^k = w^k$
    \For {$i =1,...,P$}
    \State {$\tilde f^k_{i}= f_{i}(\widetilde w_{i-1}^k)$}
        \State $\tilde d_i^k = -\nabla f_i^k(\widetilde w_{i-1}^k)$
        \State $\widetilde w_{i}^k= \widetilde w_{i-1}^k + {\zeta^k} \tilde d_i^k$
    \EndFor
    \State {\bf Output} $\widetilde w^k = \widetilde w^k_P$, $d^k = \displaystyle\sum_{i=1}^P\tilde d_i^k${, $\tilde f^k=\displaystyle\sum_{i=1}^P\tilde f^k_{i}$}
\end{algorithmic}
\caption{\texttt{Inner\_Cycle}}
\label{alg:innercycle1}
\end{algorithm}

More formally, in \CMALight we avoid evaluating $f(w^k)$ by checking the watchdog rule condition on the approximated function:
\begin{equation}
    \label{eq:ftilde}
    \tilde f^k = \sum_{i=1}^P f_i(\tilde w^k_{i-1}) = f_1 (\tilde w^k_0) + f_2 (\tilde w^k_1) + \dots + f_P (\tilde w^k_{P-1}) 
\end{equation}

where the points $\tilde w^k_i$ are defined according to the iterations of the \texttt{Inner\_Cycle}, i.e., we replace the real value of $f(w^k)$ with the sum of the $f_p$ evaluated in the points produced by an IG epoch.

According to Lemma 1 and Proposition 2 of \cite{liuzzi2022convergence}, we can prove the following result, i.e., the approximated objective function converges to the real $f(w)$.
\begin{proposition}\label{limtildew_general}
Let $f_i (w) \in \mathcal{C}^1_L$ for all $i = 1,\dots, P$. Assume that $\{w^k\}$ is bounded and that $\lim_{k\to\infty}\zeta^k = 0$. Let $\{\widetilde w^k_i\}$ and $\{d^k\}$ be the sequences of points and directions produced by Algorithm \ref{alg:IG}. 

Then, for any limit point $\Bar w$ of $\{w^k\}$ a subset of indices $K$ exists such that 
    \begin{eqnarray}
        && \lim_{k\to \infty,k\in K} w^k=\Bar w,\label{assert_C}\\
        && \lim_{k\to \infty,k\in K} \zeta^k=0,\label{assert_D}\\
        && \lim_{k\to\infty,k\in K}\widetilde w^k_i =\bar w,\quad\mbox{for all $i=1,\dots, P$}\label{assert_A_0}\\
        && \lim_{k\to\infty,k\in K} d_k = -\nabla f(\bar w).\label{assert_B_0}\\
        && {\lim_{k\to\infty,k\in K} \tilde f^k =  f(\bar w).\label{assert_C_0}}
    \end{eqnarray}
\end{proposition}

\begin{proof}
Recalling that \ref{assert_C},\ref{assert_D},\ref{assert_A_0},\ref{assert_B_0} have already been proved in \cite{liuzzi2022convergence}, we have, from Algorithm \ref{alg:innercycle1}:
$$ \tilde f^k = \sum_{i=1}^P f_i(\tilde w^k_{i-1}) = f_1 (\tilde w^k_0) + f_2 (\tilde w^k_1) + \dots + f_P (\tilde w^k_{P-1}) 
$$

Recalling that $f(w) = \sum_{i=1}^P f_i (w)$ and that  $f_i (w) \in \mathcal{C}^1_L$ for all $i = 1,\dots, P$,  using \ref{assert_A_0}, we have:
$$ \lim_{k\to\infty,k\in K} \tilde{f}^k =  \sum_{i=1}^P f_i (\bar w) = f( \bar w)
$$
\end{proof}

\section{Convergence properties of \CMALight}
\label{sec:cmalight}
In this section we define the algorithm scheme of \CMALight and prove convergence under the assumptions above.

The Controlled Minibatch Algorithm framework has been developed aiming at improving computational efficiency with respect to batch-methods, while preserving convergence properties. 
In the first version of \CMA, the convergence is ensured by a watchdog approach, imposing a condition of sufficient decrease of the objective function. 
When this condition holds, which is the most frequent case, the point returned by the Inner Cycle is accepted and no other operations are required.
When the condition is violated, i.e. the direction produced by the Inner Cycle is not a descent direction for $f(w)$, either the step-size is adjusted, or, in extreme cases, a procedure of derivative-free linesearch is executed.

In \CMALight we want to accept as much as possible the points generated by algorithm \ref{alg:innercycle1} so that in practice we almost never need to evaluate $f(w^k)$.
In order to minimize the number of function evaluation, when the watchdog conditions are not satisfied, we perform a modified version of the Extrapolation Derivative-Free Linesearch (EDFL) procedure presented in \cite{liuzzi2022convergence}, which we name EDFL Light. 

 \begin{algorithm}[h!]
\caption{Extrapolation Derivative-Free Linesearch (\texttt{EDFL\_Light})}
\label{alg:EDFL2_cma1}
\begin{algorithmic}[1]
\State Input $(\tilde f^k,w^k,d^k,\zeta^k;\gamma,\delta)$: $w^k\in \R^n,d^k\in \R^n,\zeta^k>0$, $\gamma \in (0,1), \delta\in(0,1)$
\State {Set $j = 0, \alpha = \zeta^k, f_j = \tilde f^k$}
\If {${\tilde f^k}> {f(w^k)}-\gamma\alpha\|d^k\|^2$}
    \label{edfl:step4}    \State Set $\alpha^k = 0$ and \textbf{return} %$\alpha^k$
\EndIf
\While {${f(w^k +(\alpha / \delta) d^k)}\leq \min\{f(w^k)-\gamma \alpha\|d^k\|^2, {f_j}\}$}
     \label{edfl:step7}   
     \State {Set $f_{j+1} = f(w^k +(\alpha / \delta) d^k)$}
     \State {Set $j = j + 1$}
     \State {Set $\alpha = \alpha / \delta$}
\EndWhile
\State Set $\alpha^k = \alpha$ and \textbf{return} 
\State Output: $\alpha^k, {f(w^k +\alpha^k d^k)}$
 \end{algorithmic}
 \end{algorithm}

It is trivial to verify that the following lemma holds.

\begin{lemma}\label{EDFL1:welldefined}
     Algorithm \ref{alg:EDFL2_cma1} is well defined, i.e. it determines in a finite number of steps a scalar $\alpha^k$ such that 
\begin{equation}\label{EDFLS:suffred}
f(w^k+\alpha^kd^k) \leq f(w^k) -\gamma\alpha^k\|d^k\|^2.
\end{equation}
\end{lemma}
\begin{proof}
First of all, we prove that Algorithm \ref{alg:EDFL2_cma1} is well-defined. In particular, we show that the while loop cannot infinitely cycle.  Indeed, if this was the case, that would imply
\begin{equation*}
    f\left(w^k +\frac{\zeta^k}{\delta^j} d^k\right) \leq f(w^k) -\gamma\left(\frac{\zeta^k}{\delta^j}\right)\|d^k\|^2  \quad \forall j=1,2,\dots
\end{equation*}
so that for $j\to \infty$ we would have $(\zeta^k/\delta^j)\to \infty$ which contradicts the compactness of ${\cal L}(w^0)$ and continuity of $f$.

In order to prove (\ref{EDFLS:suffred}), we consider separately the two cases i) $\alpha^k = 0$ and ii) $\alpha^k > 0$. In the first case, (\ref{EDFLS:suffred}) is trivially satisfied.

In case ii),  the condition at step {6}  implies
\[
  f(w^k+\zeta^kd^k) \leq f(w^k) -\gamma\zeta^k\|d^k\|^2.
\]
Then, if $\alpha^k = \zeta^k$, the condition is satisfied. Otherwise, if $\alpha^k > \zeta^k$, the stopping condition of the while loop, along with the fact that we already proved that the while loop could not infinitely cycle, is such that we have
\[
   f(w^k+\alpha^k d^k) \leq f(w^k) - \gamma\alpha^k\|d^k\|^2
\]
which again proves (\ref{EDFLS:suffred}), thus concluding the proof.
\end{proof}
\medskip

Therefore, the scheme of \CMALight is provided in \ref{alg:CMA1}.

We do not assume any further assumptions on the norm of $\nabla f(w)$ and we do not need for any analytical expression relating $f(w^k)$ to $\tilde f^k$, but, in order to prove the convergence, we impose conditions ensuring the boundness of the output sequence $\{w^k\}$. 
The underlying idea of \CMALight is avoiding the evaluation of $f(w^k)$ without losing the control on the behaviour of the output sequence.
In Step 6 we check that the estimate of the objective function is below its initial value and the same condition is imposed every time the EDFL is executed, in steps 15 and in step 17.
This additional control is necessary because we do not assume any particular relation between $f(w^k)$, depending only on the point $w^k$, and $\tilde f^k$, depending also on the gradients $\nabla f_i(\omega^k)$, as proved in Lemma 1 of \cite{liuzzi2022convergence}. 
In fact, this makes easy to prove the following.

\begin{lemma}
\label{lemma_levelsets_det}
Let $\{w^k\}$ be the sequence of points produced by algorithm \CMALight. Then, for all $k=0,1,\dots$, at least one of the following conditions is satisfied:
\begin{itemize}
    \item[-] $\tilde f^k \le f(\omega^0)$
    \item[-] $f(\omega^k) \le f(\omega^0)$
\end{itemize}
\end{lemma}
\begin{proof}
By definition of algorithm \CMALight\  and recalling Lemma \ref{EDFL1:welldefined},  for every iteration $k$, we have that one of the following cases happens:
\begin{enumerate}
    \item[i)] Step 7 is executed, i.e., we move to a new point such that $\tilde f^k \le \min \{\phi^{k-1} - \gamma \zeta^k, f(\omega^0) \}$.
    \item[ii)] Step 10 is executed. We move to a point such that $\tilde f^k \le f(\omega^0)$, or we set $\alpha^k = 0$, i.e., we do not move from $\omega^k$.
    \item[iii)] Step 14 is executed. After the linesearch, we move to a point such that either $\tilde f^k \le f(\omega^0)$ or  $f(\omega^k + \tilde \alpha^k d^k) = f(\omega^{k+1}) \le f(\omega^0)$ or we do not move if neither of the conditions is satisfied.
    \item[iv)] Step 16 is executed. After the linesearch we move to a point such that $f(\omega^k + \tilde \alpha^k d^k) = f(\omega^{k+1}) \le f(\omega^0)$ or we do not move. 
\end{enumerate}

The proof follows by induction. Let's consider $k=0$. For $k=0$ both the conditions are trivially satisfied. We prove that this is true for $k=1$.

If Step 7 is executed, we find $\omega^1$ such that $\tilde f^1 < \tilde f^0 = f(\omega^0)$. 
If Step 10 is executed we either move $\omega^1$ such that $\tilde f^1 \le f(\omega^0)$, or we have $\omega^1 = \omega^0$, which implies that $f(\omega^1) = f(\omega^0)$.
If step 14 is executed we  either accept the step-size returned by the linesearch and find $\omega^1$ such that $f(\omega^0) \le f(\omega^1)$, either we decrease step-size and accept $\omega^1$ such that $\tilde f^1 \le f(\omega^0)$, or we have again $\omega^1 = \omega^0$, which implies that $f(\omega^1) = f(\omega^0)$.
Eventually, if step 16 is executed we either accept the step-size returned by the linesearch and find $\omega^1$ such that $f(\omega^1) \le f(\omega^0)$, or again $\omega^1 = \omega^0$. So, we have that the lemma holds for $k=1$.

Let's now assume that the lemma is true for a generic $k$. Recalling the possible cases, we can easily proof it is true for $k+1$. In fact we will  find $\omega^{k+1}$ such that one of the following holds:
\begin{itemize}
    \item[a)] $f(\omega^{k+1}) \le f(\omega^0)$
    \item[b)] $\tilde f^{k+1} \le f(\omega^0)$
    \item[c)] $\omega^{k+1} = \omega^k$
\end{itemize}

If a) or b) holds, then the lemma holds too by definition. If c) holds, we have  $f(\omega^{k+1}) = f(\omega^k)$ and $\tilde f^{k+1} = \tilde f^k$, so that in both case the lemma still holds.
\end{proof}
\medskip

Now, we can prove the boundness of $\{ w^k \}$, which, differently from the \CMA case, is not simply guaranteed by the compactness of $\mathcal{L} (w^0)$. 
The following lemma proves that Lemma \ref{lemma_levelsets_det} ensures that $\{ w^k \}$ cannot diverge even if for some $k$ it happens that $f(w^k) > f(w^0)$.

\begin{lemma}
\label{lemma_omegalim}
Let $\{w^k\}$ be the sequence of points produced by algorithm \CMALight. Then, $\{w^k\}$ is limited.
\end{lemma}

\begin{proof}
Recalling lemma \ref{lemma_levelsets_det}, we have that the set of index $\{0,1,\dots\}$ can be represented with the following partition:
\begin{itemize}
    \item $K_1 = \{k=0,1,\dots: \tilde f^k \le f(\omega^0), f(\omega^k) > f(\omega^0)\}$
    \item $K_2 = \{k=0,1,\dots: \tilde f^k > f(\omega^0), f(\omega^k) \le f(\omega^0)\}$
    \item $K_3 = \{k=0,1,\dots: \tilde f^k \le f(\omega^0), f(\omega^k) \le f(\omega^0)\}$
\end{itemize}

Let's now suppose by contradiction that $\{\omega^k\}$ is not limited, i.e., there exists an infinite index set $K \subseteq \{0,1,\dots\}$ exists such that
\begin{equation}\label{normtoinfty}
\displaystyle \lim_{\substack{k \rightarrow \infty \\ k \in K}} \vert \vert \omega^k \vert \vert = \infty
\end{equation}

We want to show that $K_1\cap K$ is an infinite index set. Let us assume by contradiction that this is not the case. 
If this is the case, when $k\to\infty$, $k\in K$, then it must hold that $k\in K_2$ or $k\in K_3$. This implies that for $k\in K$ and sufficiently large, $k\in K_2\cup K_3$, that is $w^k\in {\cal L}(w^0)$.  This is in contrast with \eqref{normtoinfty} and assumption \ref{ass:compact}. 

Thus, we  have that:
$$ \displaystyle \lim_{\substack{k \rightarrow \infty \\ k \in K_1 \cap K}} \vert \vert \omega^k \vert \vert = \infty
$$

Recall that $\{\tilde f^k\}_{ K_1 \cap K}$ is limited both from below and above by non-negativity and by the definition of the set $K_1$. 
From the definition of $\tilde f_k$ and the non-negativity of $f_j$ for all $j$, we can write
\begin{equation}\label{boundftilde}
\tilde{f}^k = \sum_{i=1}^P f_i(\tilde{\omega}^k_{i-1}) = f_{1}(\tilde{\omega}^k_{0}) + f_{2}(\tilde{\omega}^k_{1}) + \dots + f_{P}(\tilde{\omega}^k_{P-1}) \ge f_{1} (\tilde{\omega}^k_{0})
\end{equation}
where $\tilde{\omega}^k_{0} = \omega^k$. 
Since, by coercivity of $f_{1}$, 
$$ \lim_{\substack{k \rightarrow \infty \\ k \in K_1 \cap K}} f_{1}({\omega}^k) = \infty,
$$
an index $\bar k\in K_1\cap K$ exists such that
$f_1(w^{\bar k}) > f(w^0)$. Then, by \eqref{boundftilde}, we also have
\[
\tilde f^{\bar k} \geq f_1(w^{\bar k}) > f(w^0)
\]
which contradicts the fact that, by definition of $K_1$, $\tilde f^k\leq f(w^0)$ for all $k\in K_1 \cap K$ and concludes the proof.
\end{proof}
\medskip

Thanks to this result, we only need to prove that $\lim_{k \rightarrow \infty} \zeta^k = 0$ in order to apply Proposition \ref{limtildew_general} for the convergence proof.
For this purpose we have introduced the monotonically non-increasing sequence of values $\{\phi^k\}$, which does not affect the computational efficiency of \CMALight\ but compensates the fact the both $\{\tilde f^k\}$ and $\{f(w^k)\}$ are not assumed to be non-increasing.

\begin{proposition}\label{zetazero}
Let $\{\zeta^k\}$ be the sequence of steps produced by Algorithm \CMALight, then
\[
 \lim_{k\to\infty}\zeta^k = 0.
\]
\end{proposition}
\begin{proof}
In every iteration, either $\zeta^{k+1} = \zeta^k$ or $\zeta^{k+1} = \theta\zeta^k < \zeta^k$. Therefore, the sequence $\{\zeta^k\}$ is monotonically non-increasing. Hence, it results
\[
\lim_{k\to\infty} \zeta^k = \bar\zeta \geq 0.
\]
Let us suppose, by contradiction, that $\bar\zeta > 0$. If this were the case, there should be an iteration index $\bar k\geq 0$ such that, for all $k\geq\bar k$, $\zeta^{k+1} = \zeta^k = \bar\zeta$. 
Namely, for all iterations $k\geq \bar k$, step \ref{cma2:step11} or \ref{step18_CMA2} are always executed.

Let us now denote by
\begin{eqnarray*}
  K^\prime & = & \{k:\ k\geq\bar k,\ \mbox{and step \ref{cma2:step11} is executed}\}\\
  K^{\prime\prime} & = & \{k:\ k\geq\bar k,\ \mbox{and step \ref{step18_CMA2} is executed}\}  
\end{eqnarray*}

Let first proof that  $K^\prime$ cannot be infinite. If this was the case, we would have that infinitely many times
\[
   \gamma\bar\zeta \leq \phi^k - \tilde f^{k+1} {= \phi^k - \phi^{k+1}}
\]
since, for $k\in K'$, the instructions of the algorithm imply that $\phi^{k+1} = \tilde f^{k+1}$. Note that, by the instructions of the algorithm $\{\phi^k\}$ is a monotonically non increasing sequence. Also, for all $k$, $\phi^k\geq 0$. Hence, the sequence $\{\phi^k\}$ is bounded from below hence it is convergent to a limit $\bar \phi$. Then, every sub-sequence of $\{\phi^k\}$ is converging to the same limit. In particular, 
\[
\lim_{k\to\infty,k\in K'}\phi_k = \lim_{k\to\infty,k\in K'}\phi^{k+1} = \bar\phi.
\]
That is to say that 
\[
\lim_{k\to\infty,k\in K'}\phi^k-\phi^{k+1} = 0
\]
which is in constrast with $\bar\zeta > 0$.
\par\smallskip

Thus, an index $\hat k$ exists such that,   $k\in K''$ for $k\geq\hat k$. Thanks to Lemma \ref{EDFL1:welldefined},  for $k\geq\hat k$, we have
\begin{equation}\label{seqfw}
 f(w^{k+1}) \leq f(w^k+\alpha^kd^k) \leq {f(w^k)} -\gamma\alpha^k\|d^k\|^2,
\end{equation}
that is the sequence  $\{f(w^k)\}$ is definitely monotonically non-increasing and bounded below since $f(w) \geq 0$. Hence 
\[
\lim_{k\to\infty}f(w^k) = \lim_{k\to\infty} f(w^{k+1}) = \bar f.
\]
Then, 
taking the limit in relation \eqref{seqfw} we obtain:
\[
\lim_{k\to\infty} \alpha^k\|d^k\|^2 = 0.
\]
But then, for $k\geq \hat k$ sufficiently large, it would happen that 
\[
 \alpha^k\|d^k\|^2 \leq \tau\bar\zeta
\]
which means that Algorithm \CMALight would execute step \ref{cma2:step9} setting $\zeta^{k+1} = \theta\bar\zeta$ thus decreasing $\zeta^{k+1}$ below $\bar\zeta$. This contradicts our initial assumption and concludes the proof. 
\end{proof}
\bigskip

Eventually, we remark that all the hypothesis required by Proposition \ref{limtildew_general} are satisfied. Thus, we can prove the following.

\begin{proposition}
\label{convergence}
Let $\nabla f_p(w) \in \mathcal{C}^1_L \ \ \forall p=1,\dots,P$ and let $\{\omega^k\}$ be the sequence of points produced by Algorithm \CMALight. Then, $\{\omega^k\}$ admits limit points and (at least) one of them is stationary.
\end{proposition}
\begin{proof}
By the fact that $\{\omega^k\}$ is limited, we know that $\{\omega^k\}$ admits limit points.

Now, let us introduce the following set of iteration indices
\[
 K = \{k:\ \zeta^{k+1} = \theta\zeta^k\}
\]
and note that, since $\lim_{k\to\infty}\zeta^k = 0$, $K$ must be infinite.  Let $\bar w$ be any limit point of the sub-sequence $\{\omega^k\}_K$, then, by Proposition \ref{limtildew_general}, a subset $\bar K\subseteq K$ exists such that
\begin{subequations}\label{teo:cmaconv1}
\begin{eqnarray}
&& \label{subeqa}\lim_{k\to\infty, k\in \bar K} \omega^k = \bar w,\\
&& \label{subeqb}\lim_{k\to\infty, k\in \bar K} d^k = -\nabla f(\bar w),\\
&& \label{subeqc}\lim_{k\to\infty, k\in \bar K} \zeta^k = 0
\end{eqnarray}
\end{subequations}
Then, let us now split the set $\bar K$ into two further subsets, namely,
\begin{eqnarray*}
 K_1 & = & \{k\in \bar K:\ \|d^k\| \leq \tau\zeta^k\},\\
 K_2 & = & \{k\in \bar K\setminus K_1: \widetilde\alpha^k\|d^k\|^2 \leq \tau\zeta^k\}.
\end{eqnarray*}
Note that, $K_1$ and $K_2$ cannot be both finite.

First, let us suppose that $K_1$ is infinite. Then, the definition of $K_1$, (\ref{subeqb}) and (\ref{subeqc}) imply that
\[
 \lim_{k\to\infty,k\in K_1} \|d^k\| = \|\nabla f(\bar w)\| = 0
\] 
thus concluding the proof in this case.
\par\smallskip

Now, let us suppose that $K_2$ is infinite. In this case, we proceed by contradiction and assume that $\bar w$ is not stationary, i.e. $\|\nabla f(\bar w)\| > 0$. By the definition of $K_2$ and (\ref{subeqc}), we obtain
\[
\lim_{k\to\infty,k\in K_2} \widetilde\alpha^k\|d^k\|^2 = 0,
\]
which, recalling (\ref{subeqb}) and the fact that $\|\nabla f(\bar w)\| > 0$, yields
\begin{equation}\label{eq:alfatildezero}
\lim_{k\to\infty,k\in K_2} \widetilde\alpha^k = 0.
\end{equation}
Now, let us further partition $K_2$ into two subsets, namely
\begin{itemize}
    \item[i)] $\bar K_2 = \{k\in K_2:\ \widetilde\alpha^k = 0\}$
    \item[ii)] $\widehat K_2 = \{k\in K_2:\ \widetilde\alpha^k > 0\}$.
\end{itemize}

Let us first suppose that $\bar K_2$ is infinite. This means that, by the definition of algorithm EDFL, 
\[
 \frac{f(\omega^k) - f(\omega^k+\zeta^kd^k)}{\zeta^k} < \gamma \|d^k\|^2.
\]
Then, by the Mean-Value Theorem, a number $\hat\alpha^k\in (0,\zeta^k)$ exists such that
\[
 -\nabla f(\omega^k + \hat\alpha^kd^k)^Td^k < \gamma \|d^k\|^2.
\]
Taking the limit in the above relation, considering (\ref{subeqa}), (\ref{subeqb}) and (\ref{subeqc}), we obtain
\[
 \|\nabla f(\bar w)\|^2 \leq \gamma \|\nabla f(\bar w)\|^2
\]
which, recalling that $\gamma \in (0,1)$, gives
\[
\|\nabla f(\bar w)\| = 0,
\]
which contradicts $\|\nabla f(\bar w)\|>0$.
Now, let us suppose that $\widehat K_2$ is infinite. This means that 
\[
 \frac{f(\omega^k) - f(\omega^k+(\widetilde\alpha^k/\delta)d^k)}{\widetilde\alpha^k/\delta} < \gamma \|d^k\|^2.
\]
Then, reasoning as in the previous case, and considering (\ref{eq:alfatildezero}), we again can conclude that $\|\nabla f(\bar w)\| = 0$, which is again a contradiction and concludes the proof.
\end{proof}
\bigskip
 \begin{algorithm}[h!]
\begin{algorithmic}[1]
    \State Set $ \zeta^0 >0, \theta\in(0,1), \tau > 0$, $\gamma\in(0,1), \delta\in(0,1) $
    \State Let $w^0\in \R^n,k=0$
    \State Compute $f(w^0)$ and set $\tilde f^{0} =  f(w^0)$ and ${\phi^0 = \tilde f^0}$
    \For {$k=0,1\dots $}
        \State Compute $(\widetilde w^k, d^k, {\tilde f^{k+1}})$ = \texttt{Inner\_Cycle}($w^k,\zeta^k$)
   \If{{ ${\tilde f^{k+1}} \leq \min \{{\phi^{k}} - \gamma\zeta^k, f(\omega^0)\} $}} \label{cma:step10} 
    \State {Set $\zeta^{k+1} = \zeta^k$} \label{cma2:step11} {and $\alpha^k = \zeta^k$} and ${\phi^{k+1} = \tilde f^{k+1}}$
    \Else
    \If{$\|d^k\| \leq \tau\zeta^k$}
        \State \label{cma2:step9} Set $\zeta^{k+1}=\theta \zeta^k$
        and {$\alpha^k=\begin{cases}
         \zeta^k & \text{if }{\tilde f^{k+1}}\le f(w^0)\\
        0 & \text{otherwise}\\
        \end{cases}$}
        \State {Set $\phi^{k+1} = \phi^k$}
        
    \Else
         \State $(\widetilde\alpha^k, {\hat f^{k+1}}) =\texttt{EDFL\_Light}(w^k,d^k,\zeta^k;\gamma,\delta) $
         \If {$\widetilde\alpha^k\|d^k\|^2 \leq\tau\zeta^k$}
            \State \label{step19_CMA2}Set $\zeta^{k+1}=\theta \zeta^k$  and $\alpha^k=\begin{cases} 
          \widetilde\alpha^k & \text{if } \widetilde\alpha^k > 0 \text{ {and} } {\hat f^{k+1} \le f(w^0)} \\
          \zeta^k & \text{if } \widetilde\alpha^k=0 \text{ and }{\tilde f^{k+1}}\le f(w^0)\\
          0 & \text{otherwise}\\
        \end{cases}$ 
                   
         \Else
            \State Set $\zeta^{k+1}=\zeta^k,$ \label{step18_CMA2}
            and $\alpha^k = \begin{cases} 
          \widetilde\alpha^k & \text{if } \widetilde\alpha^k > 0 \text{ {and} } { \hat f^{k+1} \le f(w^0)} \\
          0 & \text{otherwise}\\
        \end{cases}$
         \EndIf
    \State {Set $\phi^{k+1} = \min \{\hat f^{k+1}, \tilde f^{k+1}, \phi^k\}$}
    \EndIf
    \EndIf
    \State Set $w^{k+1} = w^k + \alpha^k d^k$ 
   \EndFor
\end{algorithmic}
\caption{{Controlled Minibatch Algorithm Light \texttt{(CMALight)}}}
\label{alg:CMA1}
\end{algorithm}

\section{Computational experience}
In this section we report computational experiments we have carried out to assess CMALight performance against different algorithms.

We have firstly tested CMALight on the same experimental test-bed used in \cite{liuzzi2022convergence} against CMA and other state-of-the-art algorithms.

Furthermore, we have conducted tests on the Image Classification dataset CIFAR10 \footnote{\url{https://www.cs.toronto.edu/~kriz/cifar.html}} using state-of-the-art neural architectures.

\subsection{Testing against CMA}

We have trained four different deep networks, composed only by fully-connected layers and sigmoidal activation function, which we report in Table \ref{tab:neural_archs}.
For sake of simplicity, a neural network made up of $L$ fully connected hidden layers and $N$ neurons per layer is referred to as $[L \times N]$.
Neural architectures have been trained on 11 open-source datasets, reported in Table \ref{tab:datasets} with the number of features and the number of variables of the problem.
The training process has been carried out by minimizing the MSE \eqref{eq:MSE} loss function.
The hyper-parameters setting for CMA and \CMALight is reported in Table \ref{tab:hypersetting1}.

As a benchmark for \CMALight, we have used five different algorithms:
\begin{itemize}
    \item Monotone version of CMA \cite{liuzzi2022convergence}
    \item A self-implemented Incremental Gradient (IG) method with decreasing stepsize $\zeta^{k+1} = \zeta^k (1-\epsilon)$, where $\epsilon = 0.5$
    \item Three state-of-the-art optimizers: Adam \cite{Kingma2015AdamAM}, Adagard \cite{Duchi2010AdaptiveSM}, Adadelta \cite{Zeiler2012ADADELTAAA}, which were all set to their default hyper-parameters configuration.
\end{itemize}

CMA, \CMALight, and IG have been implemented on Python 3.9, using the automatic differentiation tools for backward propagation provied by the open-source library \textit{Pytorch 2.0} \footnote{\url{https://pytorch.org/}}.

Numerical tests have been conducted on Intel(R) Xeon(R) Gold 5218 CPU @ 2.30GHz, performing 5 multi-start runs and setting a time limit of 300 seconds for each problem, where for problem we mean a pair (dataset, architecture). 
Looking at Table \ref{tab:datasets}, we remark that problem dimension on this test-bed is up to one million of variables and 500 thousands of samples.

\CMALight performance as been evaluated in terms of:
\begin{itemize}
    \item computational efficiency, by measuring the average number per iteration of evaluations of the whole objective function $f(w^k)$, i.e., with respect to all the samples
    \item bias with respect to CMA, by measuring the acceptance rate of the tentative point provided by the \texttt{Inner\_Cycle}
    \item efficiency-quality trade-off, by drawing the performance profiles curves \cite{dolan2002benchmarking}
\end{itemize}

\begin{table}[]
    \centering
    
    \begin{tabular}{c|cccc|ccc}
&\multicolumn{4}{c|}{small networks} & \multicolumn{3}{c}{large networks} \\
&          $[1\times 50]$ & $[3\times 20]$ & $[5\times 50]$ & $[10\times 50]$ &$[12\times 300]$ & $[12\times 500]$ & $[30\times 500]$ \\[.4em] \hline
         $L$ & 1 &3&5&10&12&12&30\\ 
         $N$ &50&20&50&50&300&500&500\\
%         $n$ & $50 d +50$ & $20(d+1)$\\
    \end{tabular}
    \caption{Deep Network architectures}
    \label{tab:neural_archs}
\end{table}

\begin{table}
  \centering
    \begin{adjustbox}{width=\textwidth}
        \begin{tabular}{l|l|rrr|rrrr|rrr}
            \hline
            &  & &  &  & \multicolumn{7}{c}{$n$ in K } \\ 
            &   &&&&\multicolumn{4}{c|}{small network} & \multicolumn{3}{c}{large network} \\
            &Dataset & \multicolumn{1}{r}{$P$} &  & \multicolumn{1}{r|}{$d$} & \multicolumn{1}{l}{$[1\times 50]$} & \multicolumn{1}{l}{$[3\times 20]$} & \multicolumn{1}{l}{$[5\times 50]$} & \multicolumn{1}{l}{$[10\times 50]$} & \multicolumn{1}{|l}{$[12\times 300]$} & \multicolumn{1}{l}{$[12\times 500]$} & \multicolumn{1}{l}{$[30\times 500]$}  \bigstrut[b]\\
            \hline
            \parbox[t]{2mm}{\multirow{5}{*}{\rotatebox[origin=c]{90}{small}}} 
            &Ailerons \cite{keel} & 10312 &   & 41    & 2.1  & 1.64  & 12.1 & 24.6 & 935 & 2557& 7567 \bigstrut[t]\\
            & Bejing Pm25 \cite{UCI} & 31317 &  & 48    & 2.45  & 1.78  & 12.5 & 25 & 935 & 2558 & 7568 \\
            & Bikes Sharing \cite{UCI} & 13034 &   & 59    & 3.0  & 2.0  & 13 & 2.55 & 936 & 2559 & 7569\\
            & California \cite{keel} & 15480 &   & 9     & 0.5   & 1.0  & 10.5 & 23 & 933 & 2556 & 7556 \\
            & Mv \cite{keel}    & 30576 &  & 13    & 0.7   & 1.08  & 10.7 & 23.2 & 934 & 2556 & 7566 \bigstrut[b]\\
            \hline
            \parbox[t]{2mm}{\multirow{6}{*}{\rotatebox[origin=c]{90}{big}}}
            & BlogFeedback \cite{UCI} & 39297 &  & 281   & 14.1 & 6.44  & 24.1 & 36.6 & 948 & 2570 & 7580 \bigstrut[t]\\
            & Covtype \cite{UCI} & 435759 &  & 55    & 2.8  & 1.92  & 12.8 & 25.3 & 936 & 2558 & 7568 \\
            &Protein \cite{UCI}& 109313 &  & 75    & 3.8  & 2.32 & 13.8 & 26.3 & 937 & 2560 & 7570  \\
            &Sido \cite{sido}   & 9508  &   & 4933  & 247 & 99.5 & 257 & 269 & 1180 & 2802 & 7812 \\
            &Skin NonSkin \cite{UCI} & 183792 &  & 4     & 0.25   & 0.9   & 10.3 & 22.8 & 933 & 2555& 7566 \\
            & YearPredictionMSD \cite{UCI} & 386508 &  & 91    & 4.6  & 2.64  & 14.6 & 27.1 & 937 & 2561 & 7571 \bigstrut[b]\\
            \hline
    \end{tabular}%
    \end{adjustbox}
\caption{Dataset description  ($P=$ number of training samples; $d$= number of input features) and $n=$number of variables (in K) of the optimization problems corresponding to the network architecture $[L\times N]$}
  \label{tab:datasets}%
\end{table}%

\begin{table}[h!]
    \centering
    \begin{tabular}{c|ccccc}
       &  $\zeta^0$ & $  \tau$ & $\theta $ & $\gamma $ & $
     \delta $ \\\hline
         \CMA\ & 0.5 &$10^{-2}$ &0.75 & $10^{-2}$& 0.9 \\
         \CMALight\ & 0.5 & $10^{-2}$ & 0.75&  $10^{-2}$ & 0.9
    \end{tabular}
    \caption{Hyperparameters settings in \CMA\ and \CMALight\ }
    \label{tab:hypersetting1}
\end{table}

\begin{figure*}[h!]
    \hfill
    \subfigure[]{\includegraphics[width=0.48\linewidth]{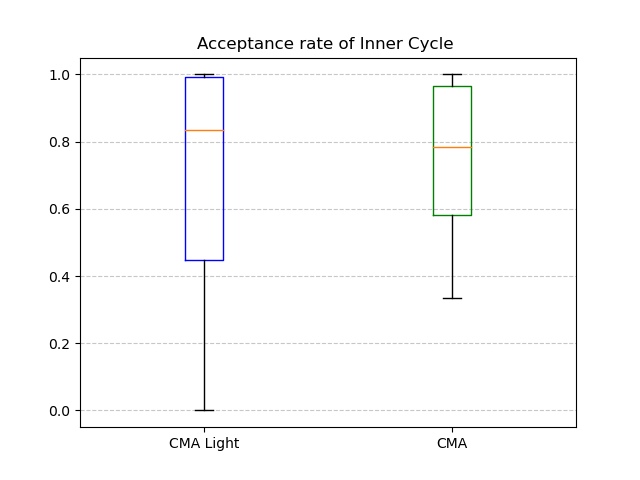}
    \label{img:ar}
    } \hfill
    \subfigure[]{\includegraphics[width=.48\linewidth]{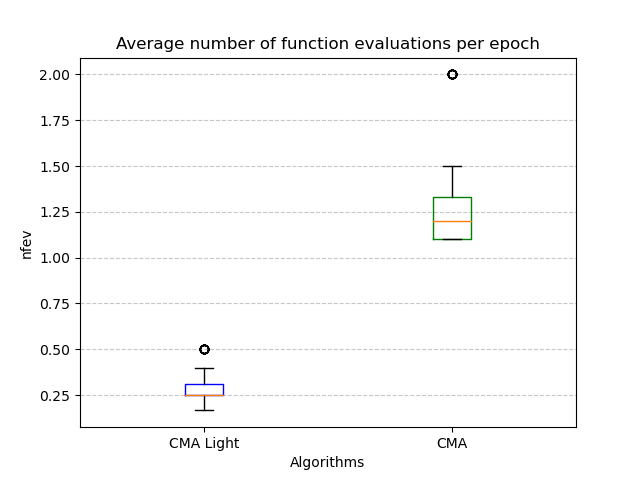}
    \label{img:nfev}
    }
    \caption{Efficiency and acceptance rate of \CMALight against CMA}
    \label{fig:FevAr}
\end{figure*}

Looking at Figure \ref{img:nfev}, we observe that, as expected, \CMALight requires significantly less average number of objective function evaluations per epoch, leading to a significant reduction of computational effort.

This holds, despite the acceptance rate of the tentative point (box-plot in Figure \ref{img:ar}) is slightly lower in \CMALight than in CMA. 
This result is not surprising, because the estimate of the objective function \eqref{eq:ftilde} is often an upper approximation for the real $f(w^k)$.
Recalling the definition, $\tilde f^k$ is computed using the sequence of points produced by Algorithm \ref{alg:innercycle1}, which means that it is equally influenced by the samples $i=1,\dots,P$.
Our computational experience shows that the anti-gradients $-\nabla f_i$ are often a descent direction even for $f(w)$, which means that the points found in the first inner iteration will correspond to worse value of the objective function.

\begin{figure*}[h!]
    \hfill
    \subfigure[]{\includegraphics[width=0.48\linewidth]{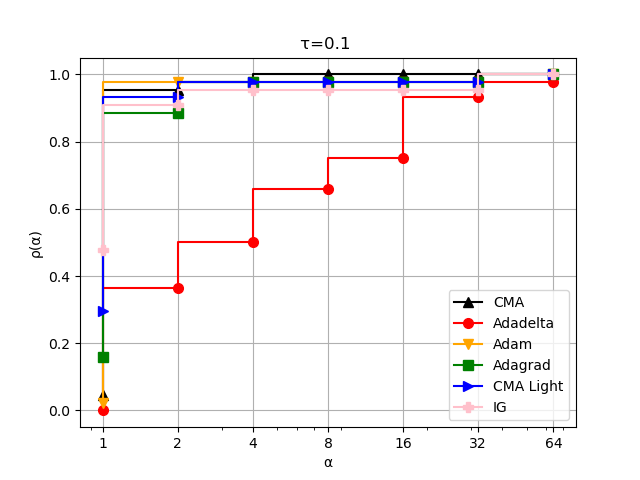}
    \label{img:PP01all}
     } \hfill
    \subfigure[]{\includegraphics[width=.48\linewidth]{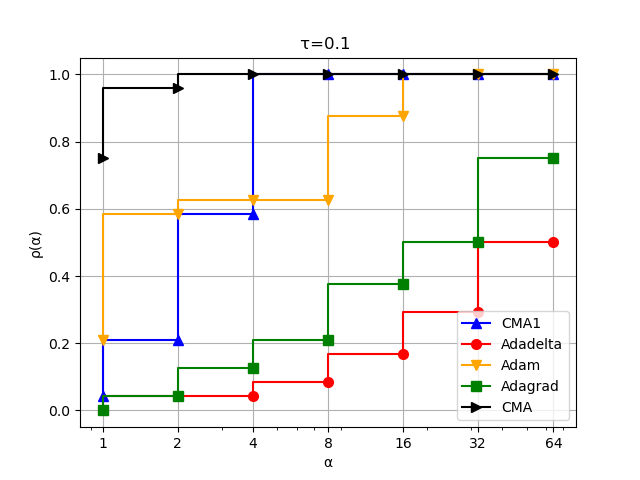}
    \label{img:PP01small}
    }
    \caption{Performance Profiles with $\tau = 0.1$ on big (a) and small (b) datasets according to classification in Table \ref{tab:datasets}}
    \label{fig:PP01}
\end{figure*}

\begin{figure*}[h!]
    \hfill
    \subfigure[]{\includegraphics[width=0.48\linewidth]{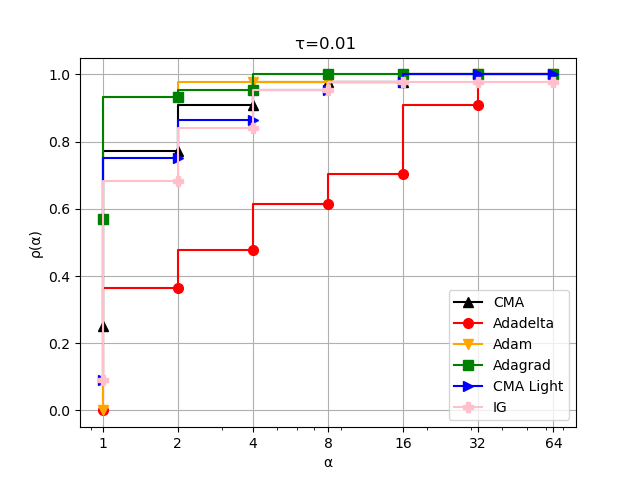}
    \label{img:PP001all}
    } \hfill
    \subfigure[]{\includegraphics[width=.48\linewidth]{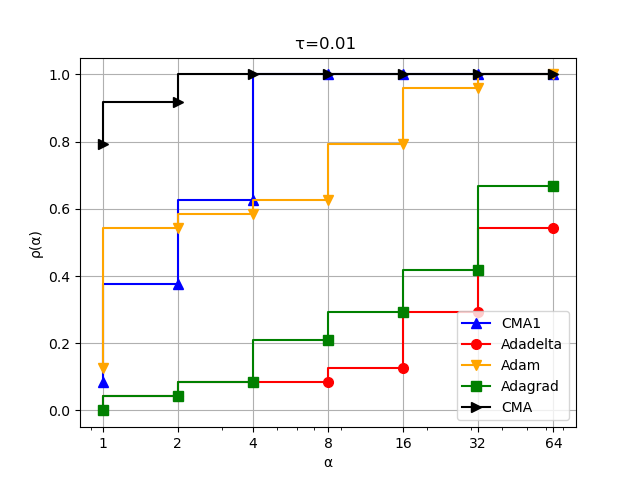}
    \label{img:PP001small}
    }
    \caption{Performance Profiles with $\tau = 0.01$ on big (a) and small (b) datasets according to classification in Table \ref{tab:datasets}}
    \label{fig:PP001}
\end{figure*}

\begin{figure*}[h!]
    \hfill
    \subfigure[]{\includegraphics[width=0.48\linewidth]{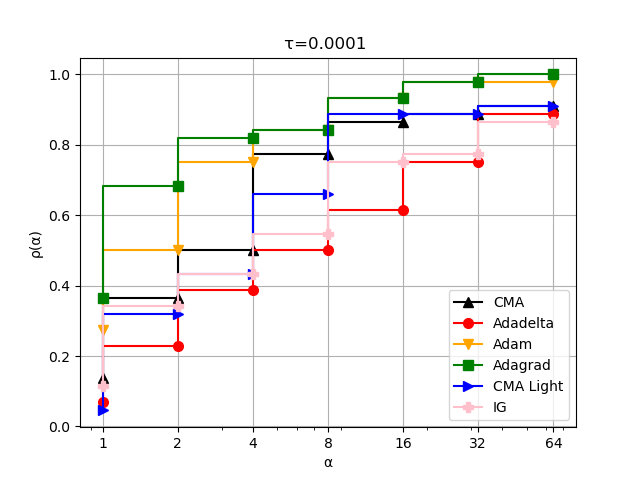}
    \label{img:PP0001all}
    } \hfill
    \subfigure[]{\includegraphics[width=.48\linewidth]{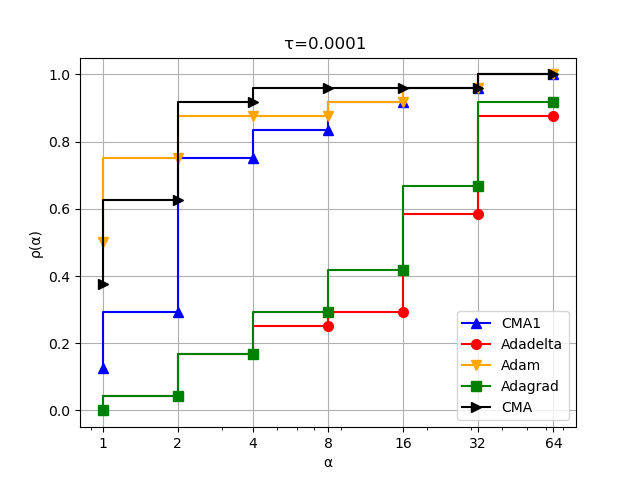}
    \label{img:PP0001small}
    }
    \caption{Performance Profiles with $\tau = 0.0001$ on big (a) and small (b) datasets according to classification in Table \ref{tab:datasets}}
    \label{fig:PP0001}
\end{figure*}

Performance profiles curves show that for low precision, in Figure \ref{fig:PP01}, \CMALight slightly outperforms Adadelta, Adagrad, IG, and CMA, when dealing with the big datasets, while on small datasets CMA plays the lion role against all the other algorithms.
When increasing the required precision, the behaviour of the algorithms on small datasets does not change significantly: CMA is still the best performing.
Nonetheless, on big datasets, the advantage of \CMALight slightly decreases in Figures \ref{fig:PP001} and \ref{fig:PP0001}.

These considerations induced us to consider, that the advantage of \CMALight could be more evident, when increasing the problem dimension toward those values, that make prohibitive to compute the objective function $f(w^k)$
\clearpage
\newpage
\subsection{Test on an Image Classification task - Early results}
Considering that \CMALight has an increasing advantage over the other algorithms on large-scale problems, we have started a new class of tests, on the dataset CIFAR10 with ResNet18 architecture \cite{he2016deep}.

CIFAR10 is an image classification dataset, made up of 50.000 $(8\times32\times32)$ images, i.e., an image has 8.192 features, belonging to ten different classes.

The target problem, in this case, is the unconstrained minimization of the cross-entropy loss (CEL \eqref{eq:CEL}) and has approximately $1.3 \times 10^{7}$ variables. 
While still considered a very basic image classification task, to the best of our knowledge, few achieved more than $94 \%$ test accuracy with ResNet18 \cite{he2016deep,SKDkaggle}, and this only using Data Augmentation \cite{van2001art}.

\begin{figure*}[h!]
    \hfill
    \subfigure[]{\includegraphics[width=0.48\linewidth]{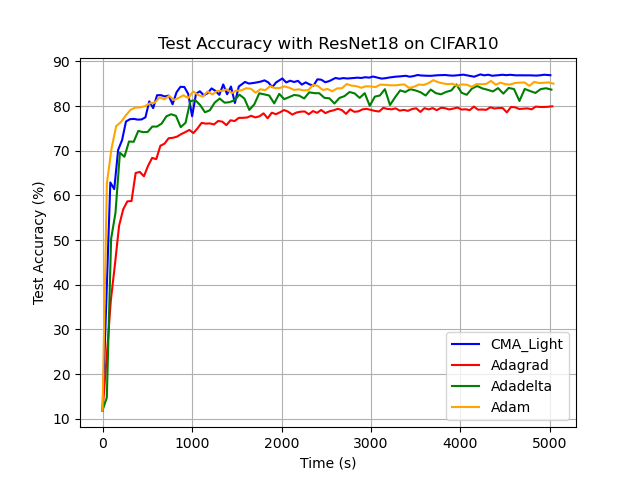}
    \label{img:resnet18acc}
    } \hfill
    \subfigure[]{\includegraphics[width=.48\linewidth]{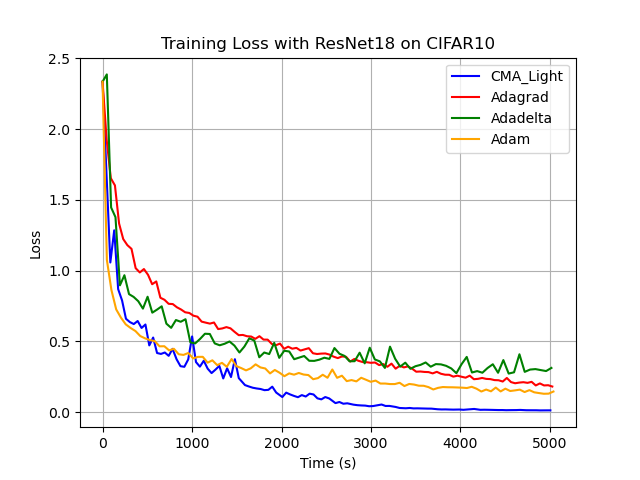}
    \label{img:resnet18loss}
    }
    \caption{Test Accuracy (a) and Training Loss (b) when training ResNet18 on CIFAR10 with different algorithms}
    \label{fig:resnet18}
\end{figure*}

Seeking for a rough idea of \CMALight performances, we performed 5 multi-start runs, training ResNet18 on CIFAR10 with Adam, Adagrad, Adadelta, and CMALight, setting a time limit of 5000 seconds and using the same device Intel(R) Xeon(R) Gold 5218 CPU @ 2.30GHz and the same hyper-parameters configuration in Table \ref{tab:hypersetting1} for \CMALight.
Results in figure \ref{fig:resnet18} shows that \CMALight has a significant advantage not only in terms of efficiency, but also in terms of quality obtained within a certain interval of time.
We remark, that the test accuracy is slightly below the state-of-the-art \cite{he2016deep,SKDkaggle}, since we did not use Data Augmentation and we used default optimizers configuration and no scheduler.

\section{Conclusions and further research}
We have developed a novel mini-batch algorithm, which converges under mild assumptions, while outperforming both CMA and other state-of-the-art solvers on large-scale datasets.
Furthermore, we have proved with an example that the algorithm can easily scale up to larger datasets and deeper neural architectures.

We aim at performing an extensive computational study on \CMALight performance on Image Classification tasks, both running tests on CIFAR10 on other deeper architectures (i.e., ResNet24, ResNet50, ResNet101, and ResNet152), and moving our research target toward harder problem, such as Image and Video Segmentation, using more complex models (e.g, Vision Transofrmers).

\clearpage
\newpage
\bibliographystyle{siamplain}
\bibliography{references}
\end{document}